\documentclass[12pt]{amsart}
\usepackage{latexsym, amssymb, amsmath}
 \usepackage{eucal}

    \newtheorem{rema}{Remark}[section]
    \newtheorem{propo}[rema]{Proposition}

   \newtheorem{theo}[rema]{Theorem}
   \newtheorem{def-theo}[rema]{Definition-Theorem}
 \newtheorem{conj}[rema]{Conjecture}
   
    \newtheorem{lemma}[rema]{Lemma}
    \newtheorem{corol}[rema]{Corollary}

	\newcommand{\nno}{\nonumber}

	\newcommand{\p}{\partial}

 \newcommand{\pf}{{\it Proof:}\hspace{2ex}}
 
 \newcommand{\epfv}{\hspace{1em}$\Box$\vspace{1em}}

\newcommand{\bC}{{\mathbb C}}
\newcommand{\bZ}{{\mathbb Z}}

\newcommand{\bQ}{{\mathbb Q}}
\newcommand{\bN}{{\mathbb N}}
\newcommand{\bT}{{\mathbb T}}

\newcommand{\BQ}{\begin{eqnarray}}
\newcommand{\EQ}{\end{eqnarray}}
\newcommand{\BQn}{\begin{eqnarray*}}
\newcommand{\EQn}{\end{eqnarray*}}
\newcommand{\BL}{\begin{align}}
\newcommand{\EL}{\end{align}}
\newcommand{\BLn}{\begin{align*}}
\newcommand{\ELn}{\end{align*}}
\newcommand{\BA}{\begin{align}}
\newcommand{\EA}{\end{align}}
\newcommand{\BAn}{\begin{align*}}
\newcommand{\EAn}{\end{align*}}
\newcommand{\wtilde}{\widetilde}

\newcommand{\Hes}{\mbox{Hes\,}}

\newcommand{\fr}{\frac}

\title[Inversion Problem and Inviscid Burgers' Equations]
{Inversion Problem, Legendre Transform and Inviscid Burgers' Equations}

    \author{Wenhua Zhao}      
    \date{\today}
    \begin{document}

\begin{abstract}
Let $F(z)=z-H(z)$ with order $o(H(z))\geq 1$ be 
a formal map from $\bC^n$ to $\bC^n$
and $G(z)$ the formal inverse map of $F(z)$.
We first study the deformation 
$F_t(z)=z-tH(z)$ of $F(z)$ and its formal inverse 
$G_t(z)=z+tN_t(z)$. (Note that $G_{t=1}(z)=G(z)$ 
when $o(H(z))\geq 2$.) 
We show that $N_t(z)$ is the unique 
power series solution of a Cauchy problem 
of a PDE,  from which we derive 
a recurrent formula for $G_t(z)$.
Secondly, 
motivated by the gradient reduction 
obtained by 
M. de Bondt, A. van den Essen
\cite{BE1}  and G. Meng \cite{M} 
for the Jacobian conjecture, 
we consider
the formal maps $F(z)=z-H(z)$ 
satisfying the gradient condition, i.e. 
$H(z)=\nabla P(z)$ for some 
$P(z)\in \bC[[z]]$ 
of order $o(P(z))\geq 2$. 
We show that, under the gradient condition,
$N_t(z)=\nabla Q_t(z)$ for some 
$Q_t(z)\in \bC[[z, t]]$ and 
the PDE satisfied by $N_t(z)$ becomes 
the $n$-dimensional inviscid Burgers' equation, 
from which a recurrent formula for $Q_t(z)$ 
also follows. Furthermore, we clarify some close 
relationships among the inversion problem, Legendre 
transform and the inviscid Burgers' equations.
In particular the Jacobian conjecture is reduced to a 
problem on the inviscid Burgers' equations.
Finally, under the gradient condition, 
we derive a binary rooted tree expansion 
inversion formula for $Q_t(z)$. 
The recurrent inversion formula and the  
binary rooted tree expansion inversion 
formula derived in this paper
can also be used as computational algorithms
for solutions of certain Cauchy problems of 
the inviscid Burgers' equations and 
Legendre transforms of 
the power series $f(z)$ 
of $o(f(z))\geq 2$.
\end{abstract}

   

 \bibliographystyle{alpha}
    \maketitle
\tableofcontents
 
\renewcommand{\theequation}{\thesection.\arabic{equation}}
\renewcommand{\therema}{\thesection.\arabic{rema}}
\setcounter{equation}{0}
\setcounter{rema}{0}
\setcounter{section}{0}

\section{\bf Introduction}\label{S1}

Let $z=(z_1, z_2, \cdots, z_n)$ and 
$F(z)=z-H(z)$  be a formal map from $\bC^n$ to 
$\bC^n$ with $o(H(z))\geq 1$
and $G(z)$ the formal inverse map of $F(z)$.
 The well-known Jacobian conjecture 
first proposed by O. H. Keller \cite{Ke} in 1939 claims that, 
{\it if $F(z)$ is a polynomial map 
with the Jacobian $j(F)(z)=1$, the inverse map $G(z)$ 
must also be a polynomial map}. Despite intense study 
from mathematicians in more than half a century, 
the conjecture is still wide open
even for the case $n=2$. For more history and known results 
on the Jacobian conjecture, 
see \cite{BCW}, \cite{E} and references there.
One of natural approaches to the Jacobian conjecture 
is to derive formulas for the inverse $G(z)$.  
In literature,  formulas which directly or indirectly 
give the formal inverse $G(z)$ 
are called inversion formulas. 
Due to many important applications in other areas, especially in 
enumerative combinatorics (See, for example, \cite{S}, \cite{Ge} 
and references there.), inversion formulas attracted 
many attentions from mathematicians much earlier 
than the Jacobian conjecture.
The first inversion formula in history 
was the Lagrange's inversion formula given 
by L. Lagrange \cite{L} in 1770, 
which provides a formula 
to calculate all coefficients of $G(z)$ 
for the one-variable case. This formula was generalized to multi-variable
cases by I. G. Good \cite{Go} in 1965. Jacobi \cite{J1} in 1830 also 
gave an inversion formula for the cases $n\leq 3$ and 
later \cite{J2} in 1833 for the general case. 
This formula now is called the Jacobi's inversion formula. 
Another inversion formula is Abhyankar-Gurjar inversion formula, 
which was first proved by Gurjar in 1974 (unpublished), 
and later Abhyankar \cite{Ab} gave a simplified proof.
By using Abhyankar-Gurjar inversion formula, 
H. Bass, E. Connell and D. Wright \cite{BCW} in 1982 and 
D. Wright \cite{Wr1} in 1989
proved the so-called Bass-Connell-Wright's tree expansion formula. 
Recently, in \cite{WZ}, this formula has been generalized
to a tree expansion formula for formal flows $F(z, t)$ generated by
$F(z)$ which provides a uniform formula for 
all the powers $F^{[m]}(z)=F(z, m)$ $(m\in \bZ)$ of $F(z)$.  
Besides the inversion formulas above, there are also many 
other inversion formulas in literature. 
See, for example, \cite{Ge}, \cite{Wr2} and references there.

Recently, M. de Bondt, A. van den Essen
\cite{BE1} and G. Meng \cite{M}
have made a breakthrough 
on the Jacobian conjecture.
They reduced the Jacobian conjecture to 
polynomial maps $F(z)=z-H(z)$ 
satisfying the {\it gradient condition}, 
i.e. $H(z)$ is the gradient $\nabla P(z)$ 
of a polynomial $P(z)$. We will refer this reduction as 
the {\it gradient reduction} and 
the condition $H(z)=\nabla P(z)$ 
the {\it gradient condition}.
One great advantage
of the gradient reduction is that, it reduces 
the Jacobian conjecture that 
involves $n$ polynomials 
to a problem that only involves a single polynomial.
Note that, by Poincar\'e lemma, a formal map $F(z)=z-H(z)$ with
$(o(H(z))\geq 1)$ satisfies 
the gradient condition if and only if the Jacobian matrix $JF(z)$ 
is symmetric. Following the terminology in \cite{BE1}, 
we also call the formal maps 
satisfying the gradient condition 
{\it symmetric} formal maps. 
For some further studies on symmetric formal maps,
see \cite{BE1}, \cite{BE2}, \cite{EW}, 
\cite{M}, \cite{Wr4}, \cite{Z1} and \cite{Z2}.

In this paper, 
we first study in Section \ref{S2} 
the deformation $F_t(z)=z-tH(z)$ of 
$F(z)=z-H(z)$ and its inverse map $G_t(z)$, 
where $t$ is a formal parameter 
which commutes with $z$.
It is easy to see that $G_t(z)$ can always 
be written as $G_t(z)=z+tN_t(z)$ 
for some $N_t(z)\in \bC[[z, t]]^{\times n}$
and $G_{t=1}(z)=G(z)$ 
when $o(H(z))\geq 2$. 
We show in Theorem \ref{T2.4} that
$N_t(z)$ is the unique solution of a Cauchy problem 
of a PDE, see Eq.\,(\ref{PDE}), (\ref{PDE-B}). 
The PDE Eq.\,(\ref{PDE})
satisfied by $N_t(z)$ has a similar form as 
the $n$-dimensional inviscid Burgers' equation 
(See \cite{R} or Eq.\,(\ref{Burgers}) 
and (\ref{Burgers2}) in this paper.). 
By solving the Cauchy problem Eq.\,(\ref{PDE}), (\ref{PDE-B}) 
recursively, we get a recurrent formula 
(See Theorem \ref{T2.7}.) 
for $N_t(z)$. This recurrent inversion formula 
not only has more computational efficiency 
in certain  situation than 
other inversion formulas, but also 
provides some new understandings on 
inversion problem. For some theoretical
consequences and applications of 
this recurrent inversion formula,  
see \cite{Wr3} and \cite{Wr4}.
Besides the main results described above,
some other properties of $N_t(z)$ 
including the one in 
Proposition \ref{P2.9} that characterizes $N_t(z)$ 
are also proved in this section. 

In Section \ref{S3}, we consider the case of 
symmetric formal maps. Let $F(z)=z-H(z)$ with 
$H(z)=\nabla (z)$ for some $P(z)\in \bC[[z]]$ 
with $o(P(z))\geq 2$. 
One can show that, in this case, 
$N_t(z)=\nabla  Q_t(z)$ for some $Q_t(z)\in \bC[[z, t]]$.
Furthermore,   
Eq.\,(\ref{PDE}) satisfied by $N_t(z)$ 
in general does become 
the $n$-dimensional inviscid Burgers' equation! 
It can be also simplified to a Cauchy problem 
Eq.\,(\ref{Cauchy-2}) 
in a single formal power series 
$Q_t(z)\in \bC[[z, t]]$ 
instead of $N_t(z)\in  
\bC[[z, t]]^{\times n}$ 
in general. 
By solving the Cauchy problem Eq.\,(\ref{Cauchy-2})
recurrently, we also get a simplified recurrent formula 
(See Proposition \ref{P3.7}.) for $Q_t(z)$. 
Some other properties of $Q_t(z)$ 
are also discussed in this section.

In Section \ref{S4}, we clarify 
some connections among the inversion problem, 
the Legendre transform and 
the inviscid Burgers' equations. 
In particular, we reduce
the Jacobian conjecture to a problem 
on the inviscid Burgers' equations, 
see Conjecture \ref{conj} and Proposition \ref{P4.2}. 
More precisely, $\nabla Q_t(z)$ is the unique power series 
solution of a Cauchy problem of the inviscid Burgers' equations
with initial condition $\nabla Q_{t=o}(z)=\nabla P(z)=H(z)$.
Note that 
the inviscid Burgers' equations are master equations for 
diffusions of airs or liquids with viscid constant $c=0$.
It is surprising for us to see that the fate of the
Jacobian conjecture is completely determined 
by behaviors of airs or liquids 
with viscid constant $c=0$.

The connection between the
inversion problem and the Legendre transform 
(See \cite{Ar} and \cite{M}) is straightforward.
For any $f(z)\in \bC[[z]]$ of order $o(f(z))\geq 2$, 
we can always write 
$f(z)=\fr 12\sum_{i=1}^n z_i^2-P(z)$ 
for some $P(z)\in \bC[[z]]$ with $o(P(z))\geq 2$.
If $\Hes (f)(0)\neq 0$,  
the Legendre transform $\bar f(z)$ of $f(z)$
is by definition
given by $\bar f(z)=\fr 12\sum_{i=1}^n z_i^2-Q(z)$, 
where $Q(z)$ is the unique formal power series with 
$o(Q(z))\geq 2$ such that the formal maps 
$F(z)=z-\nabla P(z)$ and $G(z)=z-\nabla Q(z)$ 
are inverse to each other. Hence, the
Legendre transform for formal power series $f(z)\in \bC[[z]]$ 
with $o(f(z))\geq 2$, is essentially 
the inversion problem under the gradient condition.
All results and inversion formulas derived in this paper
can also be used as computational algorithms 
for the Legendre transforms 
of formal power series 
$f(z)\in \bC[[z]]$ with $o(f(z))\geq 2$.

Finally, in Section \ref{S5}, 
by using the recurrent formula 
obtained in Proposition (\ref{P3.7}), 
we derive a binary rooted tree 
expansion inversion formula for symmetric maps, 
see Theorem \ref{T5.2}. 
Note that a tree expansion inversion formula 
for symmetric formal maps 
has been given by G. Meng \cite{M} and  
D. Wright \cite{Wr4}.
The binary rooted tree expansion inversion formula 
we derive here is different from the one 
in \cite{M} and \cite{Wr4}. 
It only involves binary rooted trees.

Two remarks are as follows. 
First, we will fix $\bC$ as our base field.
But all results, formulas as well as their proofs
given in this paper hold or work equally well 
for formal power series over any $\bQ$-algebra. 
Secondly, for convenience,
we will mainly work on the setting of 
formal power series over $\bC$. But, 
for polynomial maps or local analytic maps, 
all formal maps or 
power series involved 
in this paper are also
locally convergent.
This can be easily seen either from 
the fact that any local analytic map with non-zero Jacobian 
at the origin has a locally convergent inverse, or from 
the well-known Cauchy-Kowaleskaya theorem (See \cite{R}, for example.) 
in PDE. 

\vskip3mm

{\bf Acknowledgment}:
The author is very grateful 
to Professor David Wright 
for his encouragement and also for 
informing the author the preprint \cite{BE1} 
and some of his own recent related results. 
Great thanks also go to Professor Arno van der Essen 
for reading through the first preprint of this paper 
and pointing out many misprints. 
The author also would like to thank  
Professor Quo-Shin Chi for discussions 
on some PDE's involved in this paper and
Professor John Shareshian for informing 
the author the preprint \cite{M}.

\renewcommand{\theequation}{\thesection.\arabic{equation}}
\renewcommand{\therema}{\thesection.\arabic{rema}}
\setcounter{equation}{0}
\setcounter{rema}{0}

\section{\bf A Deformation of Formal Maps}\label{S2}

Once and for all, we fix the following notation and  conventions.

\begin{enumerate}

\item We fix $n\geq 1$ and set 
$z=(z_1, z_2, \cdots, z_n)$. For any $\bQ$-algebra $k$,
we denote by $k[z]$ (resp. $k[[z]]$) the
polynomial algebra (resp. formal power series algebra) 
over $k$ in $z_i$ $(1\leq i\leq n)$.

\item 
For any $\bQ$-algebra $k$,
by a formal map $F(z)$ from $k^n$ to $k^n$, we simply mean 
$F(z)=(F_1(z), F_2(z), \cdots, F_n(z))$
with $F_i(z)\in k[[z]]$ $(1\leq i\leq n)$. 
We denoted by $J(F)$ and $j(F)$ 
the  Jacobian matrix and the Jacobian of $F(z)$, 
respectively.

\item We denote by $\Delta$ the Laplace operator 
$\sum_{i=1}^n \frac {\p^2}{\p z_i^2}$. Note that, 
a polynomial or formal power
series $P(z)$ is said to be {\it harmonic} if $\Delta P=0$.

\item For any $k\geq 1$ and
$U(z)=(U_1(z), U_2(z), \cdots, U_k(z)) \in \bC[[z]]^{\times k}$, 
we set
\begin{align*}
o(U(z))=\min_{1\leq i\leq k} o(U_i(z))
\end{align*}
and, when $U(z)\in \bC[z]^{\times k}$,
\begin{align*}
\deg U(z)=\max_{1\leq i\leq k} \deg U_i(z).
\end{align*}
For any 
$U_t(z) \in \bC[t][z]^{\times k}
\text{ or } \bC[[z, t]]^{\times k}$ ($k\geq 1$) 
for some formal parameter $t$, the notation 
$o(U_t(z))$ and $\deg U_t(z)$
 stand for the order and 
the degree of $U_t(z)$ with respect to $z$, respectively.

\item For any $P(z)\in \bC[[z]]$, we denote by $\nabla P(z)$ the gradient 
of $P(z)$, i.e. $\nabla P=(\frac {\p P}{\p z_1}, \frac {\p P}{\p z_2}, 
\cdots, \frac {\p P}{\p z_n})$. We denote by $\Hes (P)(z)$ 
the Hessian matric of $P(z)$, 
i.e. $\Hes(P)(z)=(\frac {\p^2 P(z)}{\p z_i\p z_j})$. 

\item All $n$-vectors in this paper are supposed to be 
column vectors unless stated otherwise.
For any vector or matrix $U$, we denote by $U^\tau$ its transpose.
The standard $\bC$-bilinear form
of $n$-vectors is denoted by  $<\cdot, \cdot >$.
\end{enumerate}

\vskip6mm

In this paper, we will fix  
a formal map $F(z)$ 
from $\bC^n$ to $\bC^n$ and always assume that $F(z)$ 
has the form $F(z)=z-H(z)$ with $o(H(z))\geq 1$. Note that, 
any formal map $V: \bC^n \to \bC^n$ with $V(0)=0$ and 
$j(V)(0)\neq 0$ can be transformed into the form above by 
composing with some affine automorphisms of $\bC^n$. 

Let $t$ be a formal parameter
which commutes with $z_i$ $(1\leq i\leq n)$. We set 
$F_t(z)=z-tH(z)$. Since $F_{t=1}(z)=F(z)$, 
$F_t(z)$ can be viewed as a 
deformation of the formal map $F(z)$.
From now on, we will denote by $G(z)$ and $G_t(z)$ the formal inverses 
of $F(z)$ and $F_t(z)$, respectively. Note that,  
$G_t(z)$ can always be written as $G_t(z)=z+tN_t(z)$ for some 
$N_t(z)\in \bC[[z, t]]^{\times n}$ with $\mbox{o} (N_t(z))\geq 1$.
By uniqueness of formal inverses, we also have $G_{t=1}(z)=G(z)$. 
Furthermore, when $o(H(z))\geq 2$, 
$N_t(z)$ actually lies in $\bC[t][[z]]^{\times n}$ 
with $o(N_t(z))\geq 2$. 
This can be easily proved by using
any well-known inversion formulas, 
for example, Abhyankar-Gurjar inversion formula \cite{Ab} or
the Bass-Connell-Wright tree expansion formula \cite{BCW}. 
We will show in Theorem \ref{T2.4} that 
$N_t(z)$ is the unique solution of a Cauchy problem of PDE, 
from which we derive a recurrent formula for $N_t(z)$, 
see Theorem \ref{T2.7}.
We also discuss some other properties of $N_t(z)$ including 
the one in Proposition \ref{P2.9}, 
which characterizes $N_t(z)$, 
see Proposition \ref{P2.10}.

\begin{lemma}\label{L2.1}
For the formal power series $N_t(z)\in \bC[[z, t]]^{\times n}$ defined above, 
we have the following
identities. 
\BQ
N_t(F_t(z))&=&H(z), \label{E3.1}\\
H(G_t)&=& N_t(z). \label{E3.2}
\EQ
\end{lemma}
\pf Since $z=G_t(F_t)$, we have
\BQn
 z &=& F_t(z)+tN_t(F_t(z)), \\
  z &=& z-tH(z)+tN_t(F_t(z)).
\EQn
Therefore, 
\BQn
H(z)=N_t(F_t(z)),
\EQn
which is Eq.\,(\ref{E3.1}). 
By composing the both sides of Eq.\,(\ref{E3.1}) 
with $G_t(z)$ from right, 
we get Eq.\,(\ref{E3.2}).
\epfv

\begin{lemma}\label{L2.2}
The following statements are equivalent.
\begin{enumerate}
\item $JH(z)$ is nilpotent.
\item $\mbox{Tr } JN_t(z)=0$.
\item $JN_t(z)$ is nilpotent.
\end{enumerate}
\end{lemma}
\pf First,  by the fact $JG_t(F_t(z))=JF_t^{-1}(z)$, we have
\begin{align}
I+tJN_t(F_t)&=(I-tJH)^{-1},\nno \\
tJN_t(F_t)= I-(I-&tJH)^{-1}=tJH(I-tJH)^{-1},\nno \\
JN_t(F_t)=JH(I-&tJH)^{-1}=\sum_{k=1}^\infty JH^k(z)t^{k-1}.\label{MainEQ}
\end{align}
Therefore we have 
\BQ\label{TrEq}
\text{Tr\,} JN_t(F_t)=\sum_{k=1}^\infty  \text{Tr\,} (JH)^k t^{k-1}
\EQ
and, for any $m\geq 0$, 
\BQ\label{m-Eq}
JN_t^m(F_t)=JH^m(I-tJH)^{-m},
\EQ
since the matrices $JH$ and $(I-tJH)^{-1}$ commute with each other.

By using the fact that 
$F_t(z)$ is an automorphism of the power series algebra $\bC[[t]][[z]]$, 
we see that, 
$(1)\Leftrightarrow (2)$ follows from  Eq.\,(\ref{TrEq}) and  the fact 
that a matrix $B$ is nilpotent if and only if 
$\text{Tr\,} (B^k)=0$ $(k\geq 1)$; 
while $(1)\Leftrightarrow (3)$ follows form Eq.\,(\ref{m-Eq})
and the fact $I-tJH(z)$ is invertible in $M_n(\bC[[t]][[z]])$.
\epfv 

By Eq.\,(\ref{m-Eq}) and the fact $N_{t=1}(z)=N(z)$, 
it is easy to see that we have the following corollary.

\begin{corol}\label{C2.3}
Let $F(z)=z-H(z)$ with $o(H(z))\geq 1$ and $G(z)=z+N(z)$ with $o(N(z))\geq 1$
the formal inverse of $F(z)$. Then, for any $m\geq 1$, we have $JH^m(z)=0$ 
if and only if $JN^m(z)=0$. In particular, 
$JH(z)$ is nilpotent if and only if $JN(z)$ is.
\end{corol}

\begin{theo}\label{T2.4}
For any $H(z)\in \bC[[z]]^{\times n}$ and $N_t(z)\in \bC [t][[z]]^{\times n}$ with
$\mbox{o} (H(z))\geq 1$ and $\mbox{o} (N_t(z))\geq 1$, respectively. 
The following statements are equivalent.
\begin{enumerate}
\item The formal map $G_t(z)=z+tN_t(z)$ is the formal inverse of
$F_t(z)=z-tH(z)$.

\item $N_t(z)$ is the unique power series solution of 
the following Cauchy problem of PDE's.
\BQ
 &{}& \frac {\p N_t}{\p t}=JN_t\cdot N_t,\label{PDE}\\
 &{}& N_{t=0}(z)=H(z), \label{PDE-B}
\EQ
where $JN_t$ is the Jacobian matrix of $N_t(z)$ with respect to $z$. 
\end{enumerate}
\end{theo}

\pf First, we show $(1)\Rightarrow (2)$. By applying 
$\frac {\p}{\p t}$ to the both sides of Eq.\,(\ref{E3.1}), 
we get
\BQn
0&=&\frac{\p N_t (F_t)}{\p t}\\
&=&\frac{\p N_t}{\p t}(F_t)+
JN_t (F_t)\frac{\p F_t}{\p t}\\
&=& \frac{\p N_t}{\p t}(F_t)-
JN_t (F_t)H. 
\EQn
Therefore,
\BQn
 \frac{\p N_t}{\p t}(F_t)=
JN_t (F_t)H.
\EQn
By composing with $G_t(z)$ from right, we get
\BQn
\frac{\p N_t}{\p t}=
JN_t\cdot H(G_t)=JN_t N_t.
\EQn

Note that $G_{t=0}(z)=z$, for it is the formal inverse of 
$F_{t=0}(z)=z$.  Eq.\,(\ref{PDE-B}) 
follows immediately from Eq.\,(\ref{E3.2}) by setting $t=0$.

To show $(2)\Rightarrow (1)$, we assume that the formal 
inverse of $F_t(z)=z-tH(z)$ is given by 
$G_t(z)=z+t\widetilde N_t(z)$. By the fact proved above, we know that 
$\widetilde N_t(z)$ also satisfies Eq.\,(\ref{PDE}) 
and (\ref{PDE-B}). We will show in Proposition \ref{P2.5} below
that the power series 
solutions of the Cauchy problem  
Eq.\,(\ref{PDE}) 
and (\ref{PDE-B}) 
are actually unique. By this fact it is easy to see that
$(2)\Rightarrow (1)$ also holds.
\epfv

We define the sequence $\{ N_{[m]}(z) | m\geq 0\}$ by writing 
\begin{align}\label{Def-Nm}
N_t(z)= \sum_{m=1}^\infty  N_{[m]}(z)t^{m-1}.
\end{align}

\begin{propo} \label{P2.5}
Let $N_t(z)= \sum_{m=1}^\infty  N_{[m]}(z)t^{m-1}$
be a power series solution of  Eq.\,$(\ref{PDE})$ and $(\ref{PDE-B})$. Then
\BQ
N_{[1]}(z) &=& H(z), \label{N1}\\
N_{[m]}(z) &=& \frac 1{m-1} \sum_{\substack{k+l=m\\ k, l\geq 1}} 
JN_{[k]}(z)\cdot N_{[l]}(z)  \label{Nm}
\EQ
for any $m\geq 2$.
\end{propo}

\pf First, 
Eq.\,(\ref{N1}) follows immediately from Eq.\,$(\ref{PDE-B})$. 
Secondly, by Eq.\,(\ref{PDE}), we have

\BQn
\sum_{m=1}^\infty (m-1) N_{[m]}(z)t^{m-2} =
\left (\sum_{k=1}^\infty  JN_{[k]}(z) t^{k-1}\right )
\left (\sum_{l=1}^\infty  N_{[l]}(z)t^{l-1} \right ).
\EQn
Comparing the coefficients of $t^{m-2}$ of the both sides 
of the equation above, we have
\BQn
(m-1) N_{[m]}(z)&=& \sum_{\substack{k+l=m\\ k, l\geq 1}}  
JN_{[k]}(z) \cdot N_{[l]}(z)
\EQn
for any $m\geq 2$. Hence we get Eq.\,(\ref{Nm}). 
\epfv

By using Eq.\,(\ref{N1}), (\ref{Nm}) and the mathematical induction, 
it is easy to show the following lemma.

\begin{lemma}\label{L2.6}
$(a)$ $o(N_{[m]}(z))\geq m+1$ for any $m\geq 0$.

$(b)$ Suppose $H(z)\in \bC[z]^{\times n}$, then , for any $m\geq 1$, 
 $N_{[m]}\in \bC[z]^{\times n}$ with $\deg N_{[m]}(z) \leq (\deg H-1)m+1$.

$(c)$ If $H(z)$ is homogeneous  of degree $d$, then,
 $N_{[m]}(z)$ is homogeneous of degree
$(d-1)m+1$ for any $m\geq 1$.
\end{lemma}

Note that, by Lemma \ref{L2.6}, $(a)$, the infinite sum 
$\sum_{m=1}^\infty  N_{[m]}(z)t_0^{m-1} $ makes sense for any
complex number $t=t_0$. 
In particular, when $t=1$, $G_{t=1}(z)$ gives us 
the formal inverse $G(z)$ of $F(z)$.

\begin{theo} \label{T2.7} {\bf (Recurrent Inversion Formula)} 

Let  $\{N_{[m]}(z) | m\geq 1 \}$ be the sequence defined 
by Eq.\,$(\ref{N1})$ and $(\ref{Nm})$ recursively. Then   
the formal inverse of $F(z)=z-H(z)$ is given by 
\BQ\label{1-Inv}
G(z)=z+\sum_{m=1}^\infty N_{[m]}(z).
\EQ
\end{theo}

One interesting property of $N_t(z)$ is the following proposition.
It basically says that $\{N_t(z)| t\in \bC \}$ gives 
a family of formal maps from $\bC^n$ to $\bC^n$, 
which are ``closed'' under 
the inverse operation.

\begin{propo}\label{P2.8}
For any $s\in \bC$, 
 the formal inverse of $U_{s, t}(z)=z-sN_t(z)$
is given by
$V_{s, t}(z)=z+ s N_{t+s}(z)$.
Actually, $U_{s, t}(z)= F_{t+s}\circ G_t(z)$ 
and $V_{s, t}(z)= F_{t}\circ G_{s+t}(z)$.
\end{propo}
\pf 
\BQn
F_{t+s}\circ G_t (z)&=& G_t(z)-(t+s)H(G_t(z))\\
&=& z+tN_t(z)-(t+s)N_t(z)\\
&=& z-sN_t(z)\\
&=&U_{s, t}(z).
\EQn
 Similarly, we can prove $V_{s, t}(z)= F_{t}\circ G_{s+t}(z)$.
\epfv

Another special property of $N_t(z)$ is given 
by the following proposition.

\begin{propo}\label{P2.9}
For any $U(z)\in \bC[[z]]$, the unique power series 
solution $U_t(z)$ in $z$ and $t$ of the Cauchy problem
\begin{align}\label{GPDE}
\begin{cases}
&\frac {\p U_t}{\p t} = <\nabla U_t,  N_t>, \\
&U_{t=0}(z) = U(z).
\end{cases}
\end{align} 
is given by $U_t(z)=U(z+tN_t(z))$.
\end{propo}
\pf By similar arguments as the proof of Proposition \ref{P2.5}, 
it is easy to see that the power series solution 
in  $z$ and $t$ of the Cauchy problem Eq.\,(\ref{GPDE}) is unique. 
So it will be enough to 
show that $U_t(z)=U(z+tN_t(z))$ 
is a solution of Eq.\,(\ref{GPDE}).

\begin{align*}
\frac {\p U_t}{\p t} &=
\frac {\p}{\p t}U(z+tN_t)\\
&= 
<\nabla U (z+tN_t),\,\, \frac {\p}{\p t}(z+tN_t)>\\
&= <\nabla U (z+tN_t), \,\, N_t+t\frac {\p N_t}{\p t}> \\
\intertext{Applying Eq.\,(\ref{PDE}):}
&= <\nabla U (z+tN_t), \,\, N_t+tJN_tN_t>\\
&= <\nabla U (z+tN_t), \,\, (\text{I}+tJN_t) N_t >\\
&= <(\text{I}+tJN_t)^\tau \nabla U (z+tN_t),\, N_t >\\
&= <\nabla (U(z+tN_t)), \,\, N_t>\\
&= <\nabla U_t,\, N_t>.
\end{align*}

\epfv

Actually, $N_t(z)$ is characterized by the property 
in Proposition \ref{P2.9}. 

\begin{propo}\label{P2.10}
For any $N_t(z)\in \bC [t][[z]]^{\times n}$ with $o(N_t(z))\geq 1$, the following are equivalent.
\begin{enumerate}
\item $z+tN_t(z)$ is the formal inverse of $z-tH(z)$ for some $H(z)\in \bC [[z]]^{\times n}$.
\item Proposition \ref{P2.9} holds for $N_t(z)$.
\end{enumerate}
\end{propo}

\pf First, $(1)\Rightarrow (2)$ follows from Proposition \ref{P2.9}. 
To show $(2)\Rightarrow (1)$,   
let $U_{t, i}(z)$ $(1\leq i\leq n)$ be the unique 
power series
solution of the Cauchy problem (\ref{GPDE})
with $U(z)=z_i$ and set 
$\wtilde U_t(z)=(U_{t, 1}(z), U_{t, 2}(z), \cdots, U_{t, n}(z))$.
Note that Eq.\,(\ref{GPDE}) for $U_{t, i}(z)$ $(1\leq i\leq n)$ can be written as
\begin{align}
\frac {\p \wtilde U_t}{\p t} & = J\wtilde U_t \cdot  N_t. \label{GPDE-2}
\end{align}

By Proposition \ref{P2.9}, we have
\begin{align}\label{GPDE-3}
\wtilde U_{t}(z)&=z + tN_{t}(z).
\end{align}

By applying $\fr {\p}{\p t}$ to the equation above, we get
\BQ
\fr {\p \wtilde U_t}{\p t}=N_t+t \fr {\p N_t}{\p t}.
\EQ
By combining the equation above with Eq.\,(\ref{GPDE-2}) and (\ref{GPDE-3}), 
we have 
\begin{align*}
N_t+t \fr {\p N_t}{\p t}= J\wtilde U_t \cdot  N_t=(I+tJN_t)\cdot N_t.
\end{align*}
Therefore, we have
\begin{align}
\fr {\p N_t}{\p t}=JN_t \cdot N_t.
\end{align}
Set $H(z)=N_{t=0}(z)$. By Theorem \ref{T2.4}, we see that
$(1)$ holds.
\epfv

\renewcommand{\theequation}{\thesection.\arabic{equation}}
\renewcommand{\therema}{\thesection.\arabic{rema}}
\setcounter{equation}{0}
\setcounter{rema}{0}

\section{\bf The Case of Symmetric Formal Maps} \label{S3}

Let $F(z)=z-H(z)$ with $o(H(z))\geq 1$ 
be a formal map from $\bC^n$ to $\bC^n$. 
We say  that $F(z)$ is a {\it symmetric} formal map if 
its Jacobian matrix $J(F)$ is symmetric. 
Note that, by Poincar\'e lemma, it is easy to see that 
$F(z)$ is symmetric 
if and only if it satisfies the {\it gradient condition}, 
i.e. $H(z)=\nabla P(z)$ for some $P(z)\in \bC[[z]]$.

In this section, we study the deformation $F_t(z)$ and 
its inverse map $G_t(z)$ for symmetric formal maps $F(z)$. 
Besides some new properties of $N_t(z)$,  
the main results and formulas for $N_t(z)$ obtained 
in the previous section will also be  simplified. 

We first give a different proof for the following lemma which was first   
proved in \cite{M}.

\begin{lemma}\label{L3.1}
Let $F(z)=z-H(z)$ with $o(H(z))\geq 1$ be 
a formal map with formal inverse $G(z)=z+N(z)$. Then, 
$F(z)$ is symmetric if and only if 
$G(z)$ is.
\end{lemma}

\pf We first assume that $H(z)=\nabla P(z)$ for some $P(z)\in \bC[[z]]$. 
Note that $JH(z)=\Hes (P(z))$ is symmetric.
By Eq.\,(\ref{MainEQ}),  
we see that $JN_t (F_t)$ is symmetric. 
Hence so are $JN_t(z)$ and $JN (z)=JN_{t=1}(z)$. 
By Poincar\'e Lemma, we know that 
$N(z)$ must be the gradient of some $Q(z)\in \bC[[z]]$, i.e. 
$N(z)=\nabla Q(z)$. 

By switching $H(z)$ and $N(z)$, we see that the converse also holds.
\epfv

Now, for any 
$P(z)\in \bC[[z]]$ with $o(P(z))\geq 2$, 
we consider the deformation $F_t(z)=z-t\nabla P(z)$ and 
its inverse $G_t(z)=z+tN_t(z)$.
By Lemma $\ref{L3.1}$, we know that $N_t(z)=\nabla Q_t(z)$ 
for some $Q_t(z)\in \bC[[z, t]]$ with $o(Q_t(z))\geq 2$. 
we will fix the notation $Q_t(z)$ as above 
through the rest of this paper unless stated otherwise.  

\begin{propo}\label{P3.2}
For any $P(z)\in \bC[[z]]$ with $o(P(z))\geq 2$, 
the following  are equivalent.
\begin{enumerate}
\item $\Hes P(z)$ is nilpotent.
\item $Q_t(z)$ is harmonic, i.e. $\Delta Q_t(z)=0$.
\item $\Hes Q_t(z)$ is nilpotent.
\end{enumerate}
\end{propo}

\pf $(1)\Leftrightarrow (2)$ follows from Lemma \ref{L2.2}. 
$(1)\Leftrightarrow (3)$ 
follows from Corollary \ref{C2.3}. 
Hence we also have $(2)\Leftrightarrow (3)$. 
\epfv

\begin{lemma}\label{L3.3}
Let $P(z)$, $F_t(z)$, $G_t(z)$ and $Q_t(z)$ as above. Then we have 
the following identities. 
\begin{align}
(\nabla Q_t)(F_t)&=\nabla P \label{p-QP}, \\
(\nabla P)(G_t)&=\nabla Q_t \label{p-PQ}.
\end{align}
\end{lemma}
\pf 
Since $H(z)=\nabla P$ and $N_t(z)=\nabla Q_t$ in our case, 
the lemma follows immediately from Lemma \ref{L2.1}. 
\epfv

\begin{lemma}\label{P3.4}
\begin{align}
Q_t(F_t)&=P-\frac t2 <\nabla P, \nabla P>, \label{QP}\\
P(G_t)&=Q_t+\frac t2 <\nabla Q_t, \nabla Q_t>. \label{PQ}
\end{align}
\end{lemma}

\pf 
For any $1\leq i\leq n$,  we consider
\begin{align*}
\frac {\p Q_t(F_t)}{\p z_i} &=
\sum_{j=1}^n \frac {\p Q_t}{\p z_j}(F_t) 
\frac {\p F_{t, j}(z)} {\p z_i}\\
\intertext{Applying Eq.\,(\ref{p-QP}) in Lemma \ref{L3.3}:} 
&=\sum_{j=1}^n \frac {\p P}{\p z_j} 
(\delta_{i, j}-t \frac {\p^2 P} {\p z_i\p z_j}) \\
&=\frac {\p P}{\p z_i}- t\sum_{j=1}^n \frac {\p P}{\p z_j} 
\frac {\p^2 P} {\p z_i\p z_j} \\
&=\frac {\p P}{\p z_i}- \frac t2 \frac {\p}{\p z_i} \sum_{j=1}^n
\frac {\p P}{\p z_j} \frac {\p P}{\p z_j} \\
&=\frac{\p}{\p z_i}( P-\frac t2 <\nabla P, \nabla P>).
\end{align*}

Hence, Eq.\,(\ref{QP}) holds. Eq.\,(\ref{PQ}) can be proved similarly 
by using Eq.\,(\ref{p-PQ}).
\epfv

\begin{lemma}\label{L3.5}
\begin{align}\label{E-4QP}
\frac{\p Q_t}{\p t}(F_t)&=\frac 12 <\nabla P, \nabla P>.
\end{align}
\end{lemma}
\pf By applying $\frac {\p}{\p t}$ to the both sides of  Eq.\,(\ref{QP}), 
we get
\begin{align*}
-\frac 12 <\nabla P, \nabla P> &=
\frac{\p Q_t}{\p t}(F_t)+\sum_{j=1}^n 
\frac{\p Q_t}{\p z_j}(F_t) \frac{\p F_{t, j} }{\p t}\\
&=\frac{\p Q_t}{\p t}(F_t)-\sum_{j=1}^n \frac{\p P}{\p z_j }
\frac{\p P }{\p z_j}\\
&=\frac{\p Q_t}{\p t}(F_t)-<\nabla P, \nabla P>.
\end{align*}
Hence, Eq.\,(\ref{E-4QP}) follows.
\epfv

Under the gradient condition, 
Theorem \ref{T2.4} becomes the following theorem.

\begin{theo}\label{T3.6}
For any $Q_t(z)\in \bC[[z, t]]$ with $o(Q_t(z))\geq 2$
and $P(z)\in \bC[[z]]$ with $o(P(z))\geq 2$, 
the following are equivalent.
\begin{enumerate}
\item $G_t(z)=z+t\nabla Q_t(z)$ is the formal inverse of 
$F_t(z)=z-t\nabla P(z)$. 
\item  $Q_t(z)$ is the unique power series solution of 
the following Cauchy problem of PDE's.
\begin{align}\label{Cauchy-2}
\begin{cases}
& \frac {\p Q_t(z)}{\p t} =\frac 12 <\nabla Q_t, \nabla Q_t>, \\
& Q_{t=0}(z)=P(z).
\end{cases}
\end{align}
\end{enumerate}
\end{theo}

Note that, $(2)\Rightarrow (1)$ follows from $(1)\Rightarrow (2)$ 
and the uniqueness of the power series solutions of 
the Cauchy problem Eq.\,(\ref{Cauchy-2}). 
(For a similar argument, see the proof of Theorem \ref{T2.4}.)
While the uniqueness of 
the power series solutions of Eq.\,(\ref{Cauchy-2})
can be proved by similar arguments as the proof of 
Proposition \ref{P2.5}. 
(Also see Proposition \ref{P3.7} below.)
So we only need show $(1)\Rightarrow (2)$, for which  
we here give two different proofs. 
   
\underline {\it First Proof}: 
First note that $J(\nabla Q_t)=\Hes (Q_t)$.
By replacing $N_t(z)$ by $\nabla Q_t(z)$ 
in Eq.\,(\ref{PDE}) and (\ref{PDE-B}), we get 
\begin{align}
&\nabla \frac {\p Q_t}{\p t} =\Hes(Q_t)\nabla Q_t, \label{E3.87} \\
&\nabla Q_{t=0}(z)=\nabla P(z). \label{E3.8}
\end{align} 
Since $o(P(z))\geq 2$ and $o(Q_t(z))\geq 2$, Eq.\,(\ref{E3.8}) implies 
$Q_{t=0}(z)=P(z)$. Furthermore, Eq.\,(\ref{E3.8}) implies that,
for any $1\leq i\leq n$, we have 
\begin{align*}
\frac{\p}{\p z_i}\frac {\p Q_t}{\p t}= \sum_{j=1}^n 
\fr{\p^2 Q_t}{\p z_i\p z_j}\fr{\p Q_t}{\p z_j}
=\frac 12 \frac{\p}{\p z_i} <\nabla Q_t, \nabla Q_t>.
\end{align*} 
Since $o(\frac {\p Q_t}{\p t})\geq 2$ and $o(<\nabla Q_t, \nabla Q_t>)\geq 2$,
hence the PDE in Eq.\,(\ref{Cauchy-2}) also holds.
\epfv

\underline {\it Second Proof}: By composing with $G_t(z)$ from right to 
 both sides of Eq.\,(\ref{E-4QP}) and applying  Eq.\,(\ref{p-PQ}), 
we have 
\begin{align*}
\frac {\p Q_t(z)}{\p t}& =\frac 12 <(\nabla P)(G_t), (\nabla P)(G_t)> \\
& =\frac 12 <\nabla Q_t, \nabla Q_t>. 
\end{align*}

The initial condition in Eq.\,(\ref{Cauchy-2}),
as proved in the first proof, 
follows from Eq.\,(\ref{PDE-B}) 
in Theorem \ref{T2.4}.
\epfv

We define a sequence of formal power series 
$\{ Q_{[m]}(z)\in \bC[[z]] |m\geq 1 \}$ by writing 
\BQ\label{def-Qm}
Q_t(z)=\sum_{m=1}^\infty  Q_{[m]}(z)t^{m-1}. 
\EQ
By replacing $Q_t(z)$ by the sum above and comparing
the coefficients of $t^{m-1}$ $(m\geq 1)$ 
in Eq.\,(\ref{Cauchy-2}), it is easy to show that 
we have the following 
recurrent formula for the 
formal power series  $\{ Q_{[m]}(z)\in \bC[[z]] |m\geq 1 \}$. 

\begin{propo}\label{P3.7}
We have the following 
recurrent formula for $Q_t(z)$. 
\begin{align}
Q_{[1]}(z)&=P(z), \label{E-Recur-1.1} \\
Q_{[m]}(z)&=\frac 1{2(m-1)} 
\sum_{\substack {k, l\geq 1 \\ k+l=m}}
<\nabla Q_{[k]}(z), \nabla Q_{[l]}(z)> \label{E-Recur-1.2}
\end{align}
for any $m\geq 2$. 
In particular, when $P(z)$ is a polynomial, 
$Q_{[m]}(z)$ $(m \geq 1)$ are also polynomials. 
\end{propo}

For a uniform non-recurrent formula for 
$Q_t^k(z)$ $(k\geq 1)$ under the condition that 
$\Hes (P)$ is nilpotent, see \cite{Z1}.

\renewcommand{\theequation}{\thesection.\arabic{equation}}
\renewcommand{\therema}{\thesection.\arabic{rema}}
\setcounter{equation}{0}
\setcounter{rema}{0}

\section{\bf Relationships with Legendre 
Transform and the Inviscid Burgers' Equations} \label{S4}

In this section, we clarify some close relationships
of the inversion problem for symmetric formal maps
with the Legendre transform and  
the inviscid Burgers' equations. In particular, 
we reduce the Jacobian conjecture to a problem on 
the Cauchy problem Eq.\,(\ref{C-Burgers}), whose PDE is 
the simplified version of the inviscid Burgers' equations 
under the gradient condition.

First let us recall Legendre transform (See \cite{M} and \cite{Ar}.). 
Let $f(z)\in \bC[[z]]$ with $o(f(z))\geq 2$ and $\Hes (f)(0)\neq 0$. 
Then the formal Legendre transform $\bar f(z)$ of $f(z)$ by definition 
is the unique $\bar f(z) \in \bC[[z]]$ with $o(\bar f(z))\geq 2$ such that
the inverse map of the formal map $\nabla f: \bC^n\to \bC^n$ 
is given by $\nabla \bar f$. 
Note that, for any $f(z)\in \bC[[z]]$ of order $o(f(z))\geq 2$, 
one can always write 
$f(z)=\fr 12\sum_{i=1}^n z_i^2-P(z)$ 
for some $P(z)\in \bC[[z]]$ with $o(P(z))\geq 2$.
If $\Hes (f)(0)\neq 0$, 
it is easy to check that
the Legendre transform $\bar f(z)$ of $f(z)$ 
is given by $\bar f(z)=\fr 12 \sum_{i=1} z_i^2 +Q(z)$
for some $\bar Q(z) \in \bC[[z]]$ with $o(Q(z))\geq 2$.
Hence the Legendre transform for $f(z)\in \bC[[z]]$ 
with $o(f(z))\geq 2$ is essentially the inversion problem 
under the gradient condition. Therefore,
the recurrent inversion formula 
in Proposition \ref{P3.7} and 
the binary rooted tree expansion formula 
in Theorem \ref{T5.2} 
that will be derived in next section 
can also be used as computational algorithms 
for the Legendre transform for formal power series 
$f(z)\in \bC[[z]]$ of $o(f(z))\geq 2$.

Next we consider some relationships of the inversion problem 
for symmetric formal maps with the inviscid Burgers' equations. 
The Burgers' equations (See \cite{R}) are 
 master equations in Diffusion theory.
Recall that the $n$-dimensional 
inviscid Burgers' equation is usually written as

\BQ\label{Burgers}
\frac {\p U_t}{\p t}(z) + (JU_t)^t(z)\cdot U_t(z)=0
\EQ
or 
\BQ\label{Burgers2}
\frac {\p U_t}{\p t}(z) = (JU_t)^\tau (z)\cdot U_t(z),
\EQ
where $U_t(z)$ is a 
$n$-vector-valued function of $(t, z)$ and 
$(JU_t) (z)$ denotes the Jacobian matrix 
of $U_t(z)$ with respect to $z$. 

Note that, for any
$n$-vector-valued function $V_t(z)$ of $(t, z)$,
$V_t(z)$ satisfies Eq.\,(\ref{Burgers}) 
if and only if $-V_t(z)$ satisfies 
Eq.\,(\ref{Burgers2}). Hence Eq.\,(\ref{Burgers}) and Eq.\,(\ref{Burgers2}) 
are equivalent to each other.
In this paper, we will refer the PDE (\ref{Burgers2}) 
as the $n$-dimensional inviscid Burgers' equation.

By comparing Eq.\,(\ref{PDE}) and Eq.\,(\ref{Burgers2}), 
we see that,
the main PDE Eq.\,(\ref{PDE}) for the general inversion problem 
without the gradient condition
is almost 
the $n$-dimensional inviscid Burgers' equation 
(\ref{Burgers2}) except the transpose part. 
More interestingly, 
under the gradient condition, we have
$JN_t(z)=\Hes (Q_t)$ which 
is symmetric and Eq.\,(\ref{PDE}) 
becomes exactly 
the $n$-dimensional inviscid Burgers' equation 
Eq.\,(\ref{Burgers2}).  
The PDE in the Cauchy problem Eq.\,(\ref{Cauchy-2})
is just a simplified version 
of the inviscid Burgers' equation 
(\ref{Burgers2}) under the assumption that
$U_t(z)=\nabla Q_t(z)$ for some function 
$Q_t(z)$ of $t$ and $z$. Motivated by the connections above, 
we formulate the following conjecture.

\begin{conj}\label{conj}
For any homogeneous polynomial $P(z)$ 
of degree $d\geq 2$ with 
the Hessian matrix $\Hes (P)$ nilpotent, 
let $U_t(z)$ be the unique power series 
solution of the following Cauchy problem of PDE's. 
\begin{align}
\begin{cases}\label{C-Burgers}
& \frac {\p U_t}{\p t}(z) = \fr 12<\nabla U_t(z), \nabla  U_t(z)>, \\
& U_{t=0}(z)=P(z).
\end{cases}
\end{align}
Then $U_t(z)$ must be a polynomial in both $z$ and $t$.
\end{conj}

\begin{propo}\label{P4.2}
Conjecture $(\ref{conj})$ above for $d=4$ 
is equivalent to the Jacobian conjecture.
\end{propo}
\pf
First, by using the gradient 
reduction in \cite{BE1} and \cite{M} and 
the homogeneous reduction in \cite{BCW} 
on the Jacobian conjecture, 
we see that the Jacobian conjecture is reduced 
to polynomial maps $F(z)=z-\nabla P(z)$ 
with $P(z)$ homogeneous of degree $d=4$. 
Secondly, since $P(z)$ is homogeneous,  
the polynomial map 
$F(z)=z-\nabla P(z)$ satisfies 
the Jacobian condition $j(F)(z)=1$ if and only if 
the Hessian matrix $\Hes (P)=J(\nabla P)$ is nilpotent. 
Then it is easy to see that 
the equivalence of Conjecture \ref{conj} and 
the Jacobian conjecture follows directly 
from Theorem \ref{T3.6}. 
\epfv
 
Since the Jacobian conjecture 
for polynomial maps $F(z)$ 
of degree $\deg F(z)\leq 2$ 
has been 
proved by S. Wang \cite{Wa},
we see that Conjecture \ref{conj} 
is true for $d=2, 3$.  
It would be very interesting 
to find some proofs 
for these results by PDE methods, 
especially for the case $d=3$. 
Understandings of Conjecture \ref{conj} for $d=3$ 
from PDE point view certainly will provide 
new insights to the Jacobian conjecture.

\renewcommand{\theequation}{\thesection.\arabic{equation}}
\renewcommand{\therema}{\thesection.\arabic{rema}}
\setcounter{equation}{0}
\setcounter{rema}{0}

\section{\bf A Binary Rooted Tree Expansion Inversion Formula}\label{S5}

In this section, we derive a binary rooted tree inversion expansion formula 
for symmetric formal maps. (See Theorem \ref{T5.2}.)
First let us fix the following notations and conventions.

By a {\it rooted tree} we mean a finite
1-connected graph with one vertex designated as 
its {\it root}.
In a rooted tree
there are natural ancestral relations between vertices.  We say a
vertex $w$ is a child of vertex $v$ if the two are connected by an
edge and $w$ lies further from the root than $v$. We define the {\it degree}
of a vertex $v$ of $T$ to be the number of its children. 
A vertex is called a {\it leaf}\/ if it has no
children.  A rooted tree $T$ 
is said to be a {\it binary rooted tree} if every non-leaf vertex of $T$ 
has exactly two children. 
When we speak of isomorphisms between rooted trees, we will
always mean root-preserving isomorphisms. 

{\bf Notation:}

Once and for all, we fix the following notation for the rest of this paper.
\begin{enumerate}
\item We let $\mathbb T$ (resp. $\mathbb B$)
be the set of isomorphism classes of all
rooted trees (resp. binary rooted trees). 
For any $m\ge1$, we let $\bT_m$ 
the set of isomorphism classes of all rooted trees 
with $m$ vertices.

\item  We call  the rooted tree with one vertex the {\it singleton}, denoted by 
$\circ$.  For convenience, we also view the empty set
 as a rooted tree, denoted by $\emptyset$.

\item For any rooted tree $T$, we set the following notation:
\begin{itemize}
\item $\text{rt}_T$ denotes the root vertex of $T$.
\item $|T|$  denotes the
 number of the vertices of $T$ and 
$l(T)$ the number of leaves.
\item $\alpha (T)$ denotes the number of the elements of the automorphism
group $\mbox{Aut}(T)$.
\item  $\widehat T$ denotes the rooted tree
obtained by deleting all the leaves of $T$. 
\end{itemize}
\end{enumerate}

For any rooted trees $T_i$ with $(i=1, 2, ..., d)$, 
we define $B_+(T_1, T_2,..., T_d)$ to 
be the rooted tree obtained by connecting all roots of  $T_i$ 
$(i=1, 2, ..., d)$ to a single new vertex, which is set to the root of  
the new rooted tree $B_+(T_1, T_2,..., T_d)$. Note that, for any 
$T_1, T_2 \in {\mathbb B}$, we have $B_+(T_1, T_2)\in {\mathbb B}$.



Next let us recall $T$-factorial $T!$ 
of rooted trees $T$, which was first introduced 
by D. Kreimer \cite{Kr}. It is defined inductively as follows. 
\begin{enumerate}
\item For the empty rooted tree $\emptyset$ and the singleton $\circ$, 
we set $\emptyset !=1$ and $\circ !=1$.
\item For any rooted tree $T=B_+(T_1, T_2, ..., T_d)$, 
we set 
\begin{align}\label{EE5.1}
T!=|T|\,T_1!\,T_2!\cdots T_d!.
\end{align}
\end{enumerate}

Note that, for the chains $C_m$ $(m\in \bN )$, i.e. 
the rooted trees with $m$ vertices and height $m-1$, we have $C_m!=m!$. 
Therefore the $T$-factorial of rooted trees can be viewed as a 
generalization of the usual factorial of natural numbers.

Now, for any binary rooted tree $T$, we set 
\begin{align}\label{E5.1}
\beta(T)= \alpha(T)\, \widehat T!.
\end{align}
\begin{lemma}\label{L5.1}
$(a)$ For any non-empty binary rooted tree $T$, we have 
\begin{align}
|T|&=2l(T)-1, \label{EE5.2}\\
|\widehat T|&=l(T)-1.\label{EE5.3}
\end{align}
$(b)$ For any $T \in \mathbb B$ with 
$T=B_+(T_1, T_2)$, we have
\begin{align}\label{E5.2}
\beta(T)=\begin{cases} 2(l(T)-1)\beta(T_1)\beta(T_2)&\text{   if $T_1\simeq T_2$,}\\
 (l(T)-1)\beta(T_1)\beta(T_2)&\text{   if $T_1\not \simeq T_2$.}
\end{cases}
\end{align}
\end{lemma}
\pf $(a)$ First note that Eq.\,(\ref{EE5.3}) follows form Eq.\,(\ref{EE5.2})
and the fact $|\widehat T|=|T|-l(T)$. Hence we only need show
Eq.\,(\ref{EE5.2}).

We use the mathematical induction
on $|T|$. When $|T|=1$, we have $T=\circ$ and 
$|T|=l(T)=1$, hence $(a)$ holds.

For any $T\in \mathbb B$ with $|T|\geq 2$. 
We write $T=B_+(T_1, T_2)$. 
Note that $T_1, T_2 \neq \emptyset$ and 
$|T_i|<|T|$ $(i=1, 2)$. By our induction assumption, we have 
\begin{align*}
|T|&=|T_1|+|T_2|+1\\
&=(2l(T_1)-1)+(2l(T_2)-1)+1\\
&=2(l(T_1)+l(T_2))-1\\
&=2l(T)-1.
\end{align*}

$(b)$ First note that, we always have 
\begin{align}\label{E5.3}
\alpha(T)=\begin{cases} 2\alpha(T_1)\alpha (T_2)&\text{   if $T_1\simeq T_2$,}\\
\alpha(T_1)\alpha(T_2)&\text{   if $T_1\not \simeq T_2$.}
\end{cases}
\end{align}
By Eq.\,(\ref{EE5.1}) and (\ref{EE5.3}), we also have 
\begin{align}\label{E5.4}
 \widehat T!=|\widehat T|\, \widehat T_1!\, \widehat T_2!=
(l(T)-1)\, \widehat T_1!\, \widehat T_2!.
\end{align}
Then, it is easy to see that Eq.\,(\ref{E5.2}) 
follows directly from 
Eq.\,(\ref{E5.1}), (\ref{E5.3}) and (\ref{E5.4}).
\epfv

Now we fix $P(z)\in \bC[[z]]$ and 
$Q_t(z)\in \bC[[z, t]]$ as in Section \ref{S3}.
We sign a formal power series $Q_T(z)\in \bC[[z]]$ 
for each non-empty binary rooted tree $T$ as follows.

\begin{enumerate}
\item For $T=\circ$, we set $Q_T(z)=P(z)$.
\item For any binary rooted tree $T=B_+(T_1, T_2)$, 
we set 
\BQn
Q_T(z)=<\nabla Q_{T_1}(z), \nabla Q_{T_2}(z)>.
\EQn
\end{enumerate}

Finally we are ready to state and prove 
the main theorem of this section.

\begin{theo} \label{T5.2}
For any $m\geq 1$, we have
\begin{align} \label{MainEq5.2}
Q_{[m]}(z)
=\sum_{\substack {T\in \mathbb B_{2m-1}}} \frac 1{\beta (T)}  Q_T(z)
=\sum_{\substack {T\in \mathbb B\\ l(T)=m} } \frac 1{\beta (T)}  Q_T(z).
\end{align}
Therefore, by Eq.\,$(\ref{def-Qm})$ we have
\begin{align}
Q_t(z)&=\sum_{T\in \mathbb B\backslash \emptyset} \frac {t^{l(T)-1}}{\beta (T)}  Q_T(z),\\
Q(z)&=\sum_{T\in \mathbb B\backslash \emptyset} \frac 1{\beta (T)}  Q_T(z).  
\end{align}
\end{theo}

\pf Note that, by Eq.\,(\ref{EE5.2}) in Lemma \ref{L5.1}, we have
\BQn
\mathbb B_{2m-1}&=&\{ T\in \mathbb B | l(T)=m\} \\
\mathbb B_{2m}&=& \emptyset,
\EQn 
for any $m\geq 1$. 
Hence the two sums in Eq.\,(\ref{MainEq5.2}) are equal to each other.

To prove Eq.\,(\ref{MainEq5.2}), we first set, for any $m\geq 1$, 
$$
V_{[m]}(z)=\sum_{\substack {T\in \mathbb B\\ l(T)=m} } \frac 1{\beta (T)}  Q_T(z)
$$
and then to show that $V_{[m]}(z)=Q_{[m]}(z)$ $(m\geq 1)$. 
By Proposition \ref{P3.7}, 
it will be enough to show that the sequence 
$\{V_{[m]}(z)\in \bC[[z]] | m\geq 1\}$ also satisfy 
Eq.\,(\ref{E-Recur-1.1}) and (\ref{E-Recur-1.2}).

For the case $m=1$, since there is only one binary rooted tree $T$ with $l(T)=1$, 
namely, $T=\circ$, we have 
$V_{[1]}(z)=Q_{T=\circ}(z)=P(z)=Q_{[1]}(z)$. Hence we have Eq.\,(\ref{E-Recur-1.1}).

For any $m\geq 2$, we consider
\begin{align*}
&{} \quad \frac 1{2(m-1)} 
\sum_{\substack {k, l\geq 1 \\ k+l=m}}
<\nabla V_{[k]}(z), \nabla V_{[l]}(z)> \\ 
&=
\sum_{\substack {T_1, T_2\in \mathbb B,\\ l(T_1)=k, l(T_2)=l,\\
k, l \geq 1,  k+l=m} } \frac 1{2(m-1)\beta (T_1) \beta (T_2)}  <\nabla Q_{T_1}(z), \nabla Q_{T_2}(z)>\\
&=
\sum_{\substack {T_1, T_2\in \mathbb B,\\ l(T_1)=k, l(T_2)=l,\\
k, l \geq 1,  k+l=m} } \frac 1{2(m-1)\beta (T_1) \beta (T_2)}   Q_{B_+(T_1, T_2)}(z) \\
\intertext{Note that, the general term 
in the sum above appears twice when $T_1\not\simeq T_2$ 
but only once when $T_1\simeq T_2$. 
By applying Eq.\,(\ref{E5.2}) in  Lemma \ref{L5.1}:}
&=
\sum_{\substack {T\in \mathbb B\\ l(T)=m} } \frac 1{\beta (T)}  Q_T(z)\\
&= V_{[m]}(z).
\end{align*}
Hence we have Eq.\,(\ref{E-Recur-1.2}).
\epfv

{\small \sc Department of Mathematics, Illinois State University, 
Normal, IL 61790-4520. }

{\em E-mail}: 
wzhao@ilstu.edu.


\begin{thebibliography}{FLM2}

\bibitem[A]{Ab} S. S. Abhyankar, {\it Lectures in algebraic geometry},
Notes by Chris Christensen, Purdue Univ., 1974.

\bibitem[Ar]{Ar} V. I. Arnord, {\it Mathematical Methods of Classical Mechanics,}
Springe-Verlag New York, Inc. 1978. [MR 0690288].
 
\bibitem[BCW]{BCW} H. Bass, E. Connell, D. Wright, {\it The Jacobian
conjecture, reduction of degree and formal expansion of the inverse}.
Bull.  Amer. Math.  Soc.  \textbf{7}, (1982), 287--330. [MR 83k:14028]. 


\bibitem[BE1]{BE1} M. de Bondt and A. van den Essen, 
{\it  A Reduction of the Jacobian Conjecture to the Symmetric Case}, Report No. 0308, 
University of Nijmegen, June, 2003. To appear in {\it Proc. of the AMS.}.

\bibitem[BE2]{BE2} M. de Bondt and A. van den Essen, 
{\it  Nilpotent Symmetric Jacobian Matrices and the Jacobian Conjecture}, Report No. 0307, 
University of Nijmegen, June, 2003. To appear in {\it J. Pure and Appl, Alg.}.


\bibitem[E]{E} A. van den Essen, {\it Polynomial automorphisms and the Jacobian conjecture}.
 Progress in Mathematics, 190. Birkhäuser Verlag, Basel, 2000. MR1790619. 


\bibitem[EW]{EW} A. van den Essen and S. Washburn, 
{\it The Jacobian Conjecture for Symmetric Jacobian Matrices}, 
J. Pure Appl. Algebra, \textbf{189} (2004), no. 1-3, 123--133. 
[MR2038568]

\bibitem[Go]{Go} I. G. Good, {\it The generalization of Lagrange's expansion 
and the enumeration of trees},  Proc. Cambridge Philos. Soc. 
 {\bf 61} (1965), 499--517. [MR 31 \#88].

\bibitem[Ge]{Ge} I. M. Gessel {\it A combinatorial proof of the
multivariable Lagrange inversion formula}, J. Combin. Theory Ser. A,
\textbf{45} (1987), 178--195. [MR 88h:05011]. 

\bibitem[J1]{J1} C. G. J. Jacobi, {\it De resolutione aequationum per
series infinitas}, J. Reine Angew. Math. \textbf{6} (1830), 257--286.

\bibitem[J2]{J2} C. G. J. Jacobi, {\it Theoria novi multiplicatoris systemati
aequationum differentialium vulgarium applicandi}, 
J. Reine Angew. Math. \textbf{27} (1844), 199-268; 
{\bf 29} (1845), 213-279, 333-376.

\bibitem[Ke]{Ke} O. H. Keller, {\it Ganze Gremona-Transformation}, 
Monats. Math. Physik {\bf 47} (1939), 299-306. 

\bibitem[Kr]{Kr} {\it Chen's iterated integral represents 
 the operator product expansion},  Adv. Theor. Math. Phys. {\bf 3} 
(1999), no. 3, 627--670. [MR 1797019]. Also hep-th/9901099. 

\bibitem[L]{L} L. de Lagrange, {\it Nouvelle m\'{e}thode pour 
r\'{e}soudre des \'{e}quations litt\'{e}rales par
le moyen des s\'{e}ries.} M\'{e}m. Acad. Roy. Sci. Belles de Berlin,
{\bf 24} (1770). 

\bibitem[M]{M} G. Meng, {\it Legendre Transform, Hessian Conjecture and Tree Formula}, 
math-ph/0308035.

\bibitem[R]{R} J. Rauch, {\it Partial Differential Equations}, 
Springer-Verlag New York Inc., 1991. [MR 1223093].


\bibitem[S]{S}  Richard P. Stanley, {\it Enumerative Combinatorics
II},
Cambridge University Press, 1999. [MR 2000k:05026]. 


\bibitem[Wa]{Wa} S. S.-S. Wang, {\it A Jacobian criterion for separability}, 
J. Algebra {\bf 65} (1980), no. 2, 453--494. [MR 83e:14010]. 

\bibitem[Wr1]{Wr1} D. Wright, {\it The tree formulas for reversion of
power series}, J. Pure Appl. Algebra, {\bf 57} (1989) 191--211. 
[MR 90d:13008]. 

\bibitem[Wr2]{Wr2} D. Wright, {\it Reversion, trees, and the Jacobian conjecture}. 
Combinatorial and computational algebra (Hong Kong, 1999), 249--267, 
Contemp. Math., 264, Amer. Math. Soc., Providence, RI, 2000. [MR 1800700].

\bibitem[Wr3]{Wr3} D. Wright, {\it The Jacobian Conjecture: Ideal Membership 
Questions and recent advances}, To appear.

\bibitem[Wr4]{Wr4} D. Wright, {\it Ideal Membership Questions Relating 
to the Jacobian Conjecture.} To appear.


\bibitem[WZ]{WZ} D. Wright and W. Zhao,
{\it D-log and formal flow for analytic isomorphisms of n-space}, 
{\it Trans. Amer. Math. Soc.}, {\bf 355}, No. {\bf 8} (2003), 3117-3141.  
[MR 1974678].         Also see math.CV/0209274.  


\bibitem[Z1]{Z1} W. Zhao, {\it Hessian Nilpotent Formal Power Series and 
Their Deformed Inversion Pairs}, math.CV/0409534. 

\bibitem[Z2]{Z2} W. Zhao, {\it Some Properties and Open Problems of 
Hessian Nilpotent Polynomials}, In preparation.



\end{thebibliography}
\end{document}